\documentstyle[12pt]{article}

\font\tenmsb=msbm10
\font\sevenmsb=msbm7
\font\fivemsb=msbm5

\newfam\msbfam
\textfont\msbfam=\tenmsb
\scriptfont\msbfam=\sevenmsb
\scriptscriptfont\msbfam=\fivemsb
\def\Bbb#1{{\fam\msbfam #1}}

\makeatletter
\ifnum\@ptsize=0  \fi \ifnum\@ptsize=2  \fi 

\newcommand\qed{{\hspace*{\fill}Q.E.D.\vskip12pt plus 1pt}} 

\newcommand\sE{{\cal E}}
\newcommand\sF{{\cal F}}
\newcommand\sG{{\cal G}}
\newcommand\sH{{\cal H}}
\newcommand\sI{{\cal I}}
\newcommand\sK{{\cal K}}
\newcommand\sL{{\cal L}}
\newcommand\sO{{\cal O}}
\newcommand\sS{{\cal S}}

\newcommand\reals{{\Bbb R}}
\newcommand\zed{{\Bbb Z}}
\newcommand\comp{{\Bbb C}}

\newcommand\pn[1]{{\Bbb P}^{#1}}
\newcommand\pnpair[2]{({\Bbb P}^{#1},\sO_{{\Bbb P}^{#1}}({#2}))} \newcommand\pnsheaf[2]{\sO_{{\Bbb P}^{#1}}({#2})} \newcommand\hirz[1]{{\Bbb F}_{#1}}

\newcommand\oline[1]{{\overline {#1}}}
\newcommand\bP{{\Bbb P}}

\newcommand\proof{{\noindent\bf Proof.\ }} 

\newtheorem{theorem}{Theorem}[section]
\newtheorem{lemma}[theorem]{Lemma}
\newtheorem{corollary}[theorem]{Corollary} \newtheorem{proposition}[theorem]{Proposition} \newtheorem{question}[theorem]{Question} \newtheorem{re}[theorem]{Remark} \newtheorem{setup}[theorem]{Setup}
\newtheorem{bo}[theorem]{Basic Observation} \newtheorem{definition}[theorem]{Definition} \newtheorem{conjecture}[theorem]{Conjecture} \newtheorem{exceptions}[theorem]{The Exceptional Cases} \newtheorem{problem}[theorem]{Problem}

\newenvironment{remark}{\begin{re}\em}{\end{re}}

\begin{document}
\title {Ample Vector Bundles and Branched Coverings} \author{Thomas Peternell\and Andrew J. Sommese} \date{Dedicated to {\em Robin Hartshorne} 
on his 60th birthday}
\maketitle

\vspace*{-0.5in}\section*{Introduction} Let $f: X \to Y$ be a finite surjective morphism of degree $d$ between $n$-dimensional connected projective manifolds $X$ and $Y$. Let $\sE := [(f_*\sO_X)/\sO_Y]^*$.
It is a beautiful theorem of Lazarsfeld \cite{La80} that when $Y$ is $\pn n$, then for any such $f$ and $X$, the
vector bundle $\sE(-1)$ is spanned. There have been many partial generalizations of this result for homogeneous $Y$, e.g., \cite{De95,De96,Ki96a,Ki96b,Ma97,KM98}. In this article we investigate spannedness
and ampleness properties for not necessarily homogeneous $Y$. 

$\sE$ has a number of tendencies towards positivity. For example, using relative duality $\sE\cong f_*\omega_{X/Y}/\sO_Y$, and thus $\sE$ is weakly positive in the sense of Viehweg, e.g., \cite{Vi95}. In particular this implies the nefness of $\sE$ for $Y$ a curve. We also
prove this for $Y$ an abelian surface.
Moreover,
$\sE\otimes\sE$ always has a nontrivial section, coming from a natural linear map from $\sE^*$ to $\sE$ induced by $f$, which is an isomorphism over a dense Zariski open set of $Y$, see Lemma (\ref{genericIso}). 

Though it follows from the weak positivity that $\sE$ is nef when $\dim Y=1$, it is easy to construct examples with $\sE$ not even nef. Indeed, if $f$
expresses $X$ as a double cover of $Y$, then $\sE\cong L$ where the discriminant locus is an element of $|2L|$. Thus if $Y$ possesses a nonnef line bundle with some even power
having a smooth section, then we can construct a covering with $\sE$ not nef. 

Even when $\sE$ is nef, ampleness can fail in a number of different ways. For example, if $Y$ possesses an unramified finite cover $f:X\to Y$ of degree $d>1$, then $\sE$ is nef with no sections, but with all Chern classes
zero (see Theorem (\ref{unramifiedCase})). Thus the best possible result for ampleness would be, ``if $\sE$ is nef, and $f$ does not factor through an unramified covering of
$Y$, then $\sE$ is ample.'' The first place to check this result is for curves, since in this case $\sE$ is always nef. It is, in fact, true for $Y=\pn 1$ by Lazarsfeld's
result, and for $Y$ an elliptic curve by a result of Debarre \cite{De96}. Our Theorem (\ref{curveTheorem}) shows that this result holds for arbitrary smooth curves.

Unfortunately, in dimensions $\ge 2$, the double cover construction shows that $\sE$ can fail to be ample, even ``if $\sE$ is nef, and $f$ does not factor through an
unbranched covering of $Y$''. Some of the restrictions imposed by the hypothesis that $\sE$ is ample for all covers are given in Theorem (\ref{rarelyAmple}) and Corollary (\ref{rarelyAmpleCor}). For example, the only smooth connected projective surfaces $S$, which stand any chance of having the property that for all finite surjective morphisms from connected projective surfaces to $S$, the associated bundle $\sE$ is ample, are $\pn 2$;
a minimal $K3$-surface with no $-2$ curves; or a surface with an ample canonical bundle, for which there are neither nontrivial unramified covers over it, nor nonfinite morphisms from it onto positive dimensional varieties.

In \S \ref{surfaceSection}, we show in Theorem (\ref{delPezzoSurfaceTheorem}) that $\sE$ is always spanned for $Y$ a smooth Del Pezzo surface with $K_Y^2\ge 5$. We show in Theorem (\ref{spannedOverHirz}) that $\sE$ is spanned over the Hirzebruch surfaces ${\Bbb F}_r$, under the added condition that if $r\ge 2$, then the section
$E$, of the projection $p: Y\to \pn 1$, with negative self-intersection does not lie in the discriminant locus of the covering. We also show in Theorem (\ref{ampleOverHirz}) that $\sE$ is ample over the Hirzebruch surfaces ${\Bbb F}_r$, under the further condition that the general fiber
$F$, of the projection $p: Y\to \pn 1$, has connected inverse image under the given covering. The conditions are also shown to be necessary.

In \S \ref{delPezzoSection}, we show in Theorem (\ref{delPezzoManifoldTheorem}) that $\sE$ is always spanned for $Y$ a smooth Del Pezzo manifold of dimension
$n\ge3$ when $b:=H^n\ge 5$. Here by definition $-K_Y=(n-1)H$ for some ample line bundle $H$ on $Y$. We also show Theorem (\ref{weakExistenceTheorem}), which states
that $h^0(\sE)\not = 0$ under very mild conditions, i.e., if $b$ times the degree of the covering is at least five. Finally we ask when $h^0(\sE\otimes (-H))\not =0$.
Using adjunction theory, we show in Theorem (\ref{minuHExistenceTheorem}) and (\ref{exceptions}) that this is so for branched covers of Del Pezzo manifolds, except in
very rare circumstances.

For our induction arguments, we were forced to consider branched covers over manifolds $Y$ with domains $X$ merely normal Cohen-Macaulay with finite irrational
locus, i.e., for each point $x$ of $X$, the local rings of germs of holomorphic (or algebraic) functions at $x$ are Cohen-Macaulay and normal, with the set of nonrational singularities of $X$ finite. One interesting aspect of our
work is that many results, e.g., Lazarsfeld's original result, hold at this level of generality.

We dedicate this paper to Robin Hartshorne. His work has significantly enriched algebraic geometry. It has also been a stimulus for our work: indeed, a main result of this paper is based primarily on Robin Hartshorne's beautiful characterization \cite{Ha71} of ample vector bundles on curves.

We thank the Department of Mathematics and the Duncan Chair of the University of Notre Dame for making our collaboration possible. We also thank Meeyoung Kim and Robert Lazarsfeld for helpful comments on this paper. In particular, Robert Lazarsfeld suggested using Castelnuovo-Mumford regularity in Theorem \ref{old1.17} to improve our original result about nefness of $\sE$ off the branch locus. He also suggested Problem \ref{lazarsfeldProblem}. The second author thanks the Department of Mathematics of the K.T.H. (Royal Institute of Technology) of Stockholm, Sweden, where this paper was finished.

\tableofcontents

\section{General results on branched coverings}In this paper we work over $\comp$. By a variety we mean a complex analytic space, which might be neither reduced or irreducible. Let $f: X\to Y$ be a surjective finite morphism between $n$-dimensional projective manifolds. Let $\sE:= [(f_*\sO_X)/\sO_Y]^*$.
Here are some of
the known results on the structure of $\sE$. \begin{enumerate}
\item If $Y = \pn n,$ then $\sE$ is ample and spanned by Lazarsfeld \cite{La80}.
\item If $Y$ is a quadric, then $\sE$ is spanned by \cite{Ki96b} and $\sE$ is ample
if $3 \leq n \leq 6$ by \cite{KM98}.
\item If $Y$ is a Grassmannian, then $\sE$ is spanned, and $k$-ample \cite{Ki96a}, and in fact ample \cite{Ma97}. \item By \cite{KM98}, it follows that if $Y$ is a Lagrangian Grassmannian, then $\sE$ is ample;
and that for some other homogeneous $Y$, $\sE$ is spanned. \item $\sE$ is ample if $Y$ is an elliptic curve and if $f$ does not factor through an \'etale map \cite{De96}. \end{enumerate} We need generalizations of Lazarsfeld's theorem for mildly singular $X$. 

The general setup of these results carries over to Cohen-Macaulay $X$. Note that reduced curves have at worst Cohen-Macaulay singularities.

\begin{setup} \label{setup}Let $f:X\to Y$ be a finite, degree $d$, morphism from a reduced
pure $n$-dimensional projective variety $X$ onto an $n$-dimensional projective manifold $Y$. Assume that $X$ has at worst Cohen-Macaulay singularities. Since $f$ is flat, $f_*(\sO_X)$ is locally free of rank $d.$ By the trace map, $f_*(\sO_X)$ has a canonical direct summand $\sO_Y$; therefore we can define a locally free sheaf $\sE_f$ of rank $d-1$ by setting $$ f_*(\sO_X) = \sE_f^* \oplus \sO_Y.$$ When the map $f$ is obvious from the context, we usually denote the sheaf $\sE_f$ as $\sE$.
\end{setup}

We make the above singularity assumption on $X$, not because we want to state the most general theorems, but rather because we will need the more general choice of $X$ in induction procedures. We will often make the added assumption that $X$ is normal with the irrational locus $\sI(X)$ finite, i.e., with
the locus of nonrational singularities finite. This assumption guarantees that the crucial vanishing theorems are still valid \cite{So85}. 

We need the following slight variant of the results in \cite{So85}. 

\begin{lemma}\label{almostKodaira}Let $X$ be a normal $n$-dimensional projective variety with at worst Cohen-Macaulay singularities. Assume that the irrational locus of $X$ is finite. Assume that $L$ is ample and $D$ is an effective (possibly empty) reduced Cartier divisor. Assume that either: \begin{enumerate}
\item $n\le 2$; or
\item $D$ is normal with finite irrational locus. \end{enumerate} Then $h^i(-L-D)=0$ for $i<n$.
\end{lemma}
\proof If $D$ is empty, this is part of the Kodaira vanishing theorem of \cite{So85}. For $n=1$, the result follows immediately from Serre duality. 

Consider the sequence:
$$0\to \omega_X\otimes L\to \omega_X\otimes L\otimes \sO_X(D)\to \omega_D\otimes L_D\to 0.$$

If $n=2$ then by Serre duality we have $h^1(\omega_D\otimes L_D)=0$ and thus using \cite{So85}, $h^i(\omega_X\otimes L\otimes \sO_X(D))=0$ for $i\ge 1$. This shows the lemma when $n=2$. 

If $n>2$ then we obtain the result by using the results of \cite{So85} on $D$ and $X$.
\qed

It is conjectured that $\sE$ is always ample if $Y$ is rational homogeneous with $b_2 = 1$ \cite{KM98} or if $Y$ is simple abelian and $f$ does not factor through an \'etale map \cite{De96}. 

First we make a technically very useful generalization of Lazarsfeld's theorem to the Cohen-Macaulay case.

\begin{theorem}\label{newLazarsfeld}
Let $f : X \to \pn n$ be a finite surjective morphism from a pure dimensional reduced projective variety $X$ with at worst Cohen-Macaulay singularities. Assume that either $n=1$ or that $X$ is normal with finite irrational locus. Then
$\sE$ is spanned. If further $X$ is either irreducible or a connected curve, then $\sE(-1)$ is spanned, and, in particular, $\sE$ is ample. \end{theorem}

\proof First assume that $n=1$. Then the spannedness of $\sE$ will follow if $h^1(\sE(-1))=0$. We have
$$h^1(\sE(-1))=h^0(\sE^*(-1))=h^0((f_*\sO_X)(-1))= h^0(f^*\pnsheaf 1{-1})=0.$$

The spannedness of $\sE(-1)$ will follow from $h^1(\sE(-2))=0$. We have $$h^1(\sE(-2)=h^0(\sE^*)=h^0(f_*\sO_X)-1=h^0(\sO_X)-1.$$ Now note that if $X$ is connected, then $h^0(\sO_X)=1$. 

To show that $\sE$ is spanned for a given $n>1$, we assume by induction that the result is true for dimensions $\le n-1$. Choose an arbitrary point $y\in \pn n$.
Choose a general $D\in |\pnsheaf n 1|$ passing through $y\in \pn n$. $\oline D:=f^{-1}(D)$ is reduced, with at worst Cohen-Macaulay singularities. By induction the sheaf $\sE_{f_{\oline D}}$ is spanned. Moreover the irrational locus is finite. Then by \cite{So85} (see Lemma (\ref{almostKodaira}))
it follows that $$H^1(\sE(-1)) =H^{n-1}(X,f^*\pnsheaf n {-n}) = 0, $$ and thus the map $H^0(\sE)\to H^0(\sE_D)\to 0$. Finally note that $\sE_{f_{\oline D}}\cong \sE_D$ (see
(\ref{simpleFact}).

To show that $\sE(-1)$ is spanned for a given $n>1$, we assume by induction that the result is true for dimensions $\le n-1$. Choose an arbitrary point $y\in \pn n$.
Choose a general $D\in |\pnsheaf n 1|$ passing through $y\in \pn n$. $\oline D:=f^{-1}(D)$ is reduced, with at worst Cohen-Macaulay singularities. Moreover \cite{So86},
$\oline D$ is irreducible if $\dim D\ge 2$ and connected if $\dim \oline D =1$. Then by \cite{So85} it follows that $$H^1(\sE(-2)) =H^{n-1}(X,f^*\pnsheaf n {-n+1}) = 0.
$$ Thus the map $H^0(\sE(-1))\to H^0(\sE_D(-1))\to 0$. Now use $\sE_{f_{\oline D}}(-1)\cong \sE_D(-1)$. \qed 

By the same argument as in Theorem (\ref{newLazarsfeld}) we obtain

\begin{theorem} Let $f : X \to Q_n$ be a finite surjective morphism from a pure dimensional reduced projective variety $X$ with at worst Cohen-Macaulay singularities onto a smooth quadric $Q_n$. Assume further that $X$ is normal with
finite irrational locus. Then
$\sE$ is spanned.
\end{theorem}

We will use (and have already used) the following easy fact several times.

\begin{lemma}\label{simpleFact} Let $f : X \to Y$ be a finite surjective morphism from a pure dimensional reduced projective variety $X$, which is locally a complete intersection, onto a smooth projective manifold $Y$. Let $C \subset Y$ be a smooth (or a locally complete intersection) curve. Then, letting $X_C := f^{-1}(C)$, $\sE \vert C = (f\vert X_C)_*(\omega_{X_C/C}) / \sO_C.$ \end{lemma} 

\proof Notice that $X_C$ might not be smooth, however $X_C \subset X$ is locally a complete intersection, so that ${\omega_{X/Y}}$ makes good sense. Then our claim follows immediately from the adjunction formulas for $C \subset Y$ and $X_C \subset X$ and from the relation for the normal bundles $ N_{X_C\vert X} = (f\vert C)^*N_{C\vert Y}$. \qed

\begin{remark} Let $\sG\sS$ (respectively $\sG\sS'$) denote the set of all projective manifolds $Y$ with the property that given a finite surjective morphism $f:X\to Y$ from a pure dimensional projective manifold $X$ (respectively normal Cohen-Macaulay variety $X$ with finite irrational locus) to $Y$, $H^0(Y,\sE_f)$ spans $\sE_f$ over a dense Zariski open set of $Y$. Then it follows that given $Y\in \sG\sS$ (respectively $\sG\sS'$), $A\times Y\in \sG\sS$ (respectively $\sG\sS'$) for $A$ either projective space or the smooth quadric. This follows by the same induction used above, and the fact that given a semiample bundle $L$ on a pure $n$-dimensional
normal irreducible Cohen-Macaulay projective variety $X$ with $\sI(X)$ finite, we have the vanishing $H^j(K_X+L)=0$ for $j>n- \kappa(L)$, where $\kappa(L)$ is the Kodaira dimension of
$L$. \end{remark}

\begin{remark} Let $f:X\to Y$ be a finite morphism from a reduced pure $n$-dimensional
Cohen-Macaulay variety $X$ onto a connected projective manifold $Y$. There are many instances when $\sE_f$ is spanned. In this case, what can we say about the $k$-ampleness
of $\sE_f$ in the sense of \cite{So78}? Clearly $\sE$ is $(n-1)$-ample if and only if $X$ is connected. To see this, just note that since $\sE_f$ is assumed spanned, it
follows that if it is not $k$-ample, then there exists an irreducible subvariety $Z\subset Y$ of dimension at least $k+1$ and a trivial summand of $\sE_{f|Z}$. Thus $h^0([f_*\sO_X]_Z)\ge 2$. So if $\sE_f$ is not $(n-1)$-ample, we have $h^0(\sO_X)=h^0(f_*\sO_X)\ge 2$, which combined with $X$ being reduced, shows that $X$ is not connected. The same argument implies the
following ``criterion'' for ampleness of $\sE_f$ when $\dim Y=2$: ``if $\sE_f$ is spanned, then $\sE_f$ is ample if the restriction of $\sE_f$ to the branch locus of $f$ is ample and if inverse images under $f$ of irreducible curves are always connected.'' The key point of the argument is that if $\dim Y=2$, then an irreducible and reduced curve $C$ on $Y$ is Cartier. Thus $C'$, the inverse image of $C$ under $f$, is Cartier. Therefore, $C'$ is Cohen-Macaulay (since effective Cartier divisors on a Cohen-Macaulay variety are Cohen-Macaulay). If $C$ is not contained in the ramification locus of $f$, it follows that $C'$ is generically reduced, and therefore (since $C'$ is Cohen-Macaulay) it is reduced. Thus we conclude that $h^0(\sE^*_{f|C})=0$ if and only if $C'$ is connected. \end{remark}

Here is a general result giving ``generic semipositivity" of $\sE$, under the assumptions in (\ref{setup}) plus the assumption that $X$ is smooth. 

\begin{bo} Let $f : X \to Y$ be a finite surjective morphism from a pure dimensional reduced projective variety $X$, which is locally a complete intersection, onto a smooth projective manifold $Y$. $f_*(\omega_{X/Y}) \simeq \sE \oplus \sO_Y.$ \end{bo} 

\proof Notice first that $f_*(\omega_{X/Y}) $ is locally free of rank $d,$ which follows immediately from the
flatness of $f$ via standard cohomology theorems. Then by relative duality $$ f_*({\omega_{X/Y}})^* \simeq f_*(\sO_X)$$ which gives the claim. \qed

Our purpose is to generally investigate how ``positive" the vector bundle $\sE$ could be. To state a general result, we recall the followingdefinition due to Viehweg.

\begin{definition} {\rm Let $Y$ be a projective manifold and $\sF$ a locally free sheaf on $Y.$ Then $\sF$ is {\it weakly positive} if and only if for all ample line bundles $H$ and all positive integers $a$ the bundle $S^{ab}\sF \otimes H^b$ is generically spanned for $b \gg 0.$ } \end{definition}

To make this definition a little more suggestive we weaken this notion in the following form

\begin{definition} {\rm We say that $\sF$ is {\it generically nef,} if there exists a countable union $A = \bigcup _i A_i$ of proper analytic sets in $X$ such that $\sF\vert C$ is nef for all curves $C \not \subset A.$} \end{definition}

For a discussion of generically nef bundles we refer to \cite{DPS99}. We clearly have the following result.

\begin{proposition} If $\sF$ is weakly positive, then $\sF$ is generically nef.\end{proposition}

\begin{theorem} Let $f:X\to Y$ be a finite surjective morphism between projective manifolds. Then the bundle $\sE$ is weakly positive, in particular it is generically nef. \end{theorem} 

\proof Having in mind (1.5), this is just a special case of a result in \cite{Vi82}. \qed

We shall see below that ``in general" $\sE$ will be neither ample nor spanned nor even nef, even if $Y$ is a surface, e.g., due to the existence of $(-1)$-curves.
Nevertheless there is some tendency towards $\sE$ having sections. In (1.22) below we see by a very simple argument that $\sE \vert C$ is nef fo the
general smooth complete intersection curve $C.$

\begin{lemma}\label{genericIso}Let $f:X\to Y$ be a finite surjective morphism between projective manifolds. Then, there is an exact sequence $$0\to \sE^*\to \sE\to \sS\to 0,$$
with the reduced support of $\sS$ equal to the image under $f$ of the ramification divisor $B\subset X$. Thus $h^0(\sE\otimes \sE) \ge 1$ and there is a section of $2\det\sE$ with zero set equal set theoretically to $f(B)$.
\end{lemma}

\proof Since the higher direct images under $f$ of coherent sheaves are $0$, the direct image of the exact sequence $$ 0\to \sO_X\to \omega_{X/Y} \to \sO_X(B)\to 0,$$ gives the desired exact sequence.

Since the map from $\sE^*$ to $\sE$
is an isomorphism over a Zariski dense open set, we conclude that we have a sheaf map from $-\det \sE$ to $\det \sE$, with cokernel having support equal to the image under $f$ of the ramification divisor $B\subset X$. Thus $2\det\sE$ has a section which vanishes on the set theoretic image $f(B)$ under $f$ of the ramification divisor $B\subset X$. \qed

There is some hope for
ampleness if ${\rm Pic}(Y) = {\Bbb Z}.$

\begin{corollary} Let $f:X\to Y$ be a finite surjective morphism between $n$-dimensional projective manifolds. Assume that $f$ is not an unramified cover. If the Picard number $\rho(Y) = 1,$ then $\det \sE$ is ample.
\end{corollary}

\proof By Lemma (\ref{genericIso}), $2\det\sE$ has a section which vanishes on the set theoretic image $f(B)$ under $f$ of the ramification divisor $B\subset X$.
Since $f$ is not an unramified cover, $f(B)$ is nonempty, and $\det \sE$ is ample. \qed

\begin{remark} Using the relative Riemann-Roch theorem it follows that given a finite degree $d$ surjective morphism between $n$-dimensional projective manifolds,
$$ c_1(\sE)=\frac{f_*(\omega_{X/Y})}{2}=\frac{f_*(B)}{2},$$ where $B$ is the branch locus of $f$. Similarly, the higher Chern classes of $\sE$ can be computed. For example,
when $n=2$,
$\displaystyle c_2(\sE)=
d\chi(\sO_Y)-\chi(\sO_X)+\frac{(\omega_Y+\det\sE)\cdot\det\sE}{2}.$ \end{remark}

Next we collect some trivial but basic results on the cohomology of $\sE$ which will be used again and again.

\begin{proposition}Let $f: X\to Y$ be a finite morphism from a pure dimensional projective variety $X$
with at worst Cohen-Macaulay singularities onto a projective manifold $Y$. Then \begin{enumerate}
\item $ h^q(X,{\omega_{X/Y}}) = h^q(Y,\sE) \oplus H^q(Y,\sO_Y).$ \item $h^0(Y,\sE) =
h^0(X,{\omega_{X/Y}}) - 1 $
\item $H^q(X,\sO_X) = H^q(Y,\sE^*) \oplus H^q(Y,\sO_Y).$ \item If $Y$ is
Fano and $X$ satisfies the additional conditions of being normal with finite irrational locus, then $H^q(Y,\sE) = 0$ for $q \geq 1.$ \item If $Y$
is Fano, and $X$ satisfies the additional conditions of being normal with finite irrational locus, and $-\omega_X = rH$ with $r \geq 2$ (so $X$ has index at least 2), then $H^q(Y,\sE(-H)) = 0$ for $q \geq 1.$ \\ \end{enumerate}
\end{proposition}

\proof (1), (2) and (3) being clear, (4) and (5) follow by Kodaira's vanishing theorem. \qed

Notice that (5) implies that if $h^q(\sO_X) > h^q(\sO_Y)$ for some $1 \leq q \leq n-1,$ then $\sE$ cannot be
Nakano positive.
We will investigate $\sE$ mostly on Fano manifolds of large index and in low dimensions. Here is
one result of more general nature.

\begin{theorem}\label{old1.17} Let $f:X\to Y$ be a finite surjective morphism from the projective manifold $X$ of dimension $n$ to the abelian variety $Y$. Then $\sE_f$ is nef.
\end{theorem}

\proof This is an easy adoption of the proof of \cite[3.21]{DPS94}. Adopting the arguments of (ibid.(ii)) word by word, we need to prove the following. \\ 

Let $G$ be a very ample line bundle on $Y.$ Then $f_*(\omega_X) \otimes G^{n+1}$ is spanned.\\

This follows immediately from noting that $f_*(\omega_X) \otimes G^{n+1}$ is Castelnuovo-Mumford $n$-regular. Indeed, let $y$ be a point of $Y$. Choose $n$
sections $s_1,\ldots,s_n$ of $G$, whose scheme-theoretic zero set includes $y$ as a nonsingular component. Let $K^\bullet$ denote the tensor the Koszul complex associated to the section $s_1\oplus\cdots\oplus s_n$ of $\displaystyle G^{\oplus n}$ with $f_*(\omega_X) \otimes
G$. Applying Kodaira vanishing to the hypercohomology of $K^\bullet$ shows that $f_*(\omega_X) \otimes
G^{n+1}$ is spanned.

Then
introduce $\lambda : Y \to Y, z \to 2z $, make the base change $\lambda^p$ and argue as in the second part of \cite{DPS94}. \\ \qed

\begin{corollary} If $Y$ is the product of an abelian variety and projective spaces, then $\sE$ is nef.
\end{corollary}

\begin{problem}[Lazarsfeld]\label{lazarsfeldProblem} Let $f:X\to Y$ be a finite surjective morphism between $n$-dimesnional projective manifolds. Is it true that $\sE_f$ is always nef modulo the discriminant locus D, i.e., $\sE \vert C$ is nef for all curves $C \not \in D$?
\end{problem}

The following summarizes the general things that we know in the unramified case. 

\begin{theorem}\label{unramifiedCase}
Let $f : X \to Y$ be a finite unramified morphism from a connected $n$-dimensional projective manifold $X$ onto a connected $n$-dimensional projective manifold $Y$. Then
\begin{enumerate}
\item $\sE\cong\sE^*$ and $h^0(\sE)=0$;
\item $\sE$ is nef;
\item the $i$-th Chern class $c_i(\sE)$ vanishes for $i\ge 1$; and \item $\sE$ has a filtration by hermitian flat vector bundles. \end{enumerate} \end{theorem}

\proof Since $\omega_X\cong f^*\omega_Y$, we have that $\omega_{X/Y}\cong \sO_X$. Taking direct images we conclude that $$
\sE^*\oplus \sO_Y\cong f_*(\sO_X)\cong f_*(\omega_{X/Y})\cong \sE\oplus \sO_Y,
$$
and therefore that $\sE^*\cong\sE$.

To see that $h^0(\sE)=0$, simply note that $$
h^0(\sE)=h^0(\sE^*)=h^0(\sO_X)-h^0(\sO_Y)=0. $$ 

To see nefness of $\sE$, it is equivalent to show nefness of $(f_*\sO_X)^*$. By definition this comes down to showing that there is no finite morphism $g:C\to Y$ of an irreducible and reduced projective curve $C$ to $Y$ with the property that there is a surjection $\displaystyle g^*[V^*]\to \sL^*\to 0$ of the pullback under $g$ of the vector bundle $V^*$ associated $(f_*\sO_X)^*$ onto $\sL^*$ a line bundle on $C$ with $\deg\sL>0$. We argue by contradiction. Assume otherwise that there is such a surjection, or equivalently that there is an injection $$
0\to \sL\to g^*(f_*\sO_X).
$$

We have an inner product defined on $f_*\sO_X$ by using the trace. Given two elements $s,t$ of the fiber of $V$ at $y\in Y$, define $(s,t)={\rm tr}(s\overline t)$,
where, letting $\sI_{f^{-1}(y)}$ denote the ideal sheaf of the fiber of $f$ over $y$, we identify $s,t$ with the corresponding elements of $\displaystyle\sO_X/\sI_{f^{-1}(y)}$.
Define the function $h:V\to \reals$ by sending $s\in V$ to $(s,s)$. This function is a plurisubharmonic exhaustion of $V$, since $f$ is unramified. In particular
the restriction of $h$ to $\sL$ is a
nonconstant plurisubharmonic function on $\sL$. But, since $\sL$ has positive degree, any plurisubharmonic function on $\sL$ is constant. This contradiction shows the nefness of $\sE$. 

Since the metric $(s,t)$ is flat, i.e., its curvature form vanishes, we conclude that the Chern classes $c_i(\sE)$ of $\sE$ with $i>0$ vanish. For the last statement we refer to [DPS94]. \qed 

\begin{corollary}
Let $f : X \to Y$ be a finite morphism from a connected $n$-dimensional projective manifold $X$ onto a connected $n$-dimensional projective manifold $Y$.
Then $f$ is unramified if and only if $c_1(\sE)=0$. \end{corollary} 

\proof By Theorem (\ref{unramifiedCase}), we can assume without loss of generality that $f$ is ramified. If $n=1$, then the result follows from the following Lemma.

\begin{lemma}\label{detEonCurve} Let $f:X\to Y$ be a finite morphism between smooth connected curves. Let $\tau$ denote the degree of the ramification divisor of $f.$ Then $f_*({\omega_{X/Y}})$ and $\sE$ are nef, and $\deg f_*({\omega_{X/Y}})= {{\tau} \over {2}};$ in particular $\det (\sE)$ is ample unless $f$ is an unramified cover.
\end{lemma}

\proof The nefness follows from (1.12). Let $c := \deg f_*({\omega_{X/Y}}) = \deg \sE.$ By Riemann-Roch we have $$ \chi (f_*({\omega_{X/Y}})) = d(1-g(Y)) + c,$$ where $g$ denotes the genus. On the other hand Riemann-Roch on $X$ gives \begin{eqnarray*}
\chi(f_*{\omega_{X/Y}}) &=& \chi({\omega_{X/Y}}) = \chi(\omega_X \otimes f^*(\omega_Y^{-1}))\\ &=& - \chi(f^*(\omega_Y^{-1})) = g(X)-1 -d(2g(Y)-2). \end{eqnarray*} Putting things together and using the
Riemann-Hurwitz formula, the claim follows. \qed

Thus we can assume that $f$ is ramified and also that $n\ge 2$. In this case
let $C$ be a smooth curve on $Y$ obtained as the intersection of $n-1$ general elements of $|L|$, where $L$ is a very ample line bundle on $Y$. By Bertini's theorem, $C$ and $X_C:=f^{-1}(C)$ are smooth. By Kodaira's vanishing theorem, $X_C$ is connected. Since $X_C$ is the intersection of ample divisors, it meets the ramification divisor nontrivially. This gives that $\deg(f_*(\omega_{X_C/C})/\sO_C)>0$. Since $f_*(\omega_{X_C/C})/\sO_C\cong \sE_C$, we conclude the contradiction that $c_1(\sE)\not=0$. \qed

\begin{proposition}\label{unf} Suppose that $f: X \to Y$ factors through an unramified cover: $f = b \circ a$ with $b: Z \to Y$ unramified of degree at least 2. Then $\sE$ cannot
be ample and moreover $\sE$ is not spanned at any point. \end{proposition}

\proof In fact, $\sE_b = (b_*(\sO_Z)/\sO_Y)^*$ is numerically flat with no sections by Theorem (\ref{unramifiedCase}). Since $\sE_b$ is a direct summand of $\sE$, the claims follow. \qed 

\begin{corollary}Let $f:X\to Y$ be a finite surjective morphism between smooth connected $n$-dimensional projective manifolds. Let $r$ denote the rank of the subsheaf of $\sE$ generated by the images of elements of $H^0(\sE)$.
If $f$ factors through an unramified cover $b:Z\to Y$ from a connected projective manifold $Z$ onto $Y$, then the sheet number of $b$ is bounded by $r+1$. In particular, if global sections of $\sE$ span $\sE$ at some point, then
\begin{enumerate}
\item $f$ does not factor
through any nontrivial unramified cover of $Y$; and \item given a smooth $C$ obtained as the intersection of general elements $A_i\in |L_i|$, for $n-1$ ample and spanned line bundles $L_1,\ldots,L_{n-1}$, we have $\sE_C$ is ample. \end{enumerate} \end{corollary}

Probably the generic spannedness of $\sE$ is unneeded in the next result. 

\begin{corollary}Let $f:X\to Y$ be a finite surjective morphism between smooth connected $n$-dimensional projective manifolds. Assume that global sections of $\sE$ span $\sE$ at some point, and let $C$ be the smooth curve obtained as the intersection of general elements $A_i\in |L_i|$, for $n-1$ ample and spanned line bundles $L_1,\ldots,L_{n-1}$. Then the fundamental group of $f^{-1}(C)$ surjects onto the fundamental group of $C$. \end{corollary} 

\proof First note that $C':=f^{-1}(C)$ is smooth by Bertini's theorem, and connected by the Kodaira vanishing theorem. Since $\sE_{C'}\cong f_{C*}(\omega_{C'/C})/\sO_C$, by Lemma (\ref{simpleFact}), $f_{C'}$ does not factor through an unramified cover of $C$. This is equivalent to having the fundamental group of $C'$ surject onto the fundamental group of $C$. \qed

\begin{remark}\label{badExamples}
We would like to give some examples that show that $\sE$, even when ample, may have very few or no sections. Let $L$ be a line bundle on a projective manifold $Y$. Assume that there exists a section $s$ of $L$ with a smooth zero set $D$. Assume that $A$ is a line bundle on $Y$ with $rA=L$. Let $f:X\to Y$ be the $r$-sheeted cover obtained by taking the $r$-th root of $s$. In this case $\displaystyle \sE
\cong \bigoplus_{j=1}^{r-1}jA$. Thus if we choose $A$ and $L$ such that $A$ is not spanned, we have an example with $\sE$ not spanned. Such examples are plentiful. For example,
\begin{enumerate}
\item We could take $L$ as the
line bundle associated to two distinct points on a curve $Y$. If $h^1(\sO_Y)>0$, then any line bundle $A$ with $2A=L$, is of degree one, and thus has at most one section. A dimension count shows that if $h^1(\sO_Y)\ge 2$, we can choose the two points so that in fact $h^0(A)=0$; \item If $Y$ is a Del Pezzo surface, i.e., if $-\omega_Y$ is ample, then it follows by a straightforward argument using Reider's theorem, that if $kA=L$ with $h^0(L)\ge 1$, $A$ not spanned, and $k\ge 2$, then $L=-k\omega_Y$, $A=-\omega_Y$, and $\omega_Y\cdot \omega_Y=1$. \item There is
a similar example for Del Pezzo threefolds. Let $Y$ denote a Del Pezzo threefold with $-\omega_Y=2H$ with $H^3=1$. Then $-\omega_Y$ is spanned, but for the double cover associated to $2H=-\omega_Y$, $\sE=H$ is spanned except at one point.
\end{enumerate}
Similar examples show that $\sE$ does not have to be nef. Let $Y$ be a Hirzebruch surface $F_r$ with $r\ge 2$, i.e., a $\pn 1$-bundle over $\pn 1$ with a section $E$ satisfying $E\cdot E=-r\le -2$. Let $f$ denote a fiber of the tautological surjection of $F_r$ to $\pn 1$. Take $L$ as the line bundle $rE+r(r-1)f$. Note that since $(r-1)[E+rf]$ is spanned we can choose a smooth divisor $D'\in |(r-1)[E+rf]|$. Since $E\cdot D'=0$, $D:=E+D'$ is a smooth divisor $\in |rE+r(r-1)f|$. Taking the $r$-sheeted cyclic cover associated to the $r$-th root of
$D$, we have
$$
\sE \cong \bigoplus_{j=1}^{r-1}jA,
$$
where $A$ the line bundle associated to $E+(r-1)f$. Note that this bundle is negative restricted to $E$, and so $\sE$ is not even nef in this case. \end{remark}

It is worth emphasizing that it is hard for $\sE$ to be ample for all branched covers over a fixed base.

\begin{theorem}\label{rarelyAmple}Let $Y$ be an $n$-dimensional connected projective manifold. Assume that for all finite surjective morphisms $f:X\to Y$ from connected projective manifolds $X$, the associated bundle $\sE$ is ample. Then:
\begin{enumerate}
\item the profinite completion of the fundamental group of $Y$ is $0$, and thus in particular the first betti number of $Y$ is $0$; and \item there exist no surjective morphisms $g: Y\to Z$ with $\dim Z\ge 1$, $Z$ projective, and $g$ having at least one positive dimensional fiber. \end{enumerate}
\end{theorem}

\proof The first assertion is clear. To see the second assertion, let $L$ be a very ample line bundle on $Z$. Choosing the double cover associated to a smooth $D\in |2g^*L|$, we have $\sE\cong g^*L$, which, by the hypothesis on $g$, cannot be ample.
\qed

\begin{corollary}\label{rarelyAmpleCor}Let $Y$ be an $n$-dimensional connected projective manifold. Assume that for all finite surjective morphisms $f:X\to Y$ from
connected projective manifolds $X$, the associated bundle $\sE$ is ample. Then either $-\omega_X$ is ample with {\rm Pic}$(Y)=\zed$ or $\omega_Y$ is nef. Moreover $\omega_Y$ is ample if $\omega_Y$ is big. \end{corollary} 

\begin{remark} Using the prediction of the Abundance Conjecture that $mK_Y$ is spanned for suitable large $m$ whenever $K_Y$ is nef, it follows in (\ref{rarelyAmpleCor}),
that either $-\omega_Y$ is ample or $\omega_Y$ is ample or $\omega_Y \equiv 0.$ In the last case $Y$ is an irreducible Calabi-Yau or symplectic manifold.
In particular, if $\dim Y = 2,$ then $Y$ is the projective plane or $\omega_Y$ is ample or $Y$ is K3 without $(-2)-$curves. If $\dim Y = 3,$ then
$Y$ is Fano with $b_2 = 1$ or $\omega_Y$ is ample or $Y$ is Calabi-Yau without any contraction.
\end{remark}

\section{Coverings of curves}\label{curveSection}In this section we prove a result generalizing the theorems of Lazarsfeld \cite{La80} for $\pn 1$ and Debarre \cite{De96} for elliptic curves. 

\begin{theorem}\label{curveTheorem}
Let $f: X\to Y$ be a finite morphism from a smooth connected projective curve $X$ onto a smooth connected projective curve $Y$. The bundle $\sE=[(f_*\sO_X)/\sO_Y]^*$ is ample if and only if $f$ does not factor through a nontrivial unramified covering of $Y$. \end{theorem} 

\begin{remark}
As is pointed out in (\ref{badExamples}), even when $\sE$ is ample, \begin{enumerate}
\item if $h^1(\sO_Y)>0$, then $\sE$ does not have to be spanned; and
\item if $h^1(\sO_Y)\ge 2$, then $h^0(\sE)$ can equal $0$. \end{enumerate} When $h^1(\sO_Y)=1$, then
$h^1(\sE)=h^0(\sE^*)=0$ and hence by the Riemann-Roch theorem, and Lemma (\ref{detEonCurve}),
$$
h^0(\sE)=\deg(\sE)=h^1(\sO_X)-1.
$$
\end{remark}

The rest of the section is devoted to the proof of Theorem (\ref{curveTheorem}).

By (1.23) we may assume that $f$ does not factor through an unramified covering. We assume that $\sE$ is not ample, and argue by contradiction. We have the following simple corollary of Hartshorne's characterization of ampleness of vector bundles on curves \cite{Ha71}. 

\begin{lemma}\label{hartshorne}
Let $V$ be a nef vector bundle on a smooth connected curve $C$. There exists a unique maximal ample vector subbundle $A$. The quotient $V/A$ has degree zero.
\end{lemma}

\proof If $V$ is ample the Lemma is trivially true. Thus we can assume that $V$ is not ample.

By Hartshorne's characterization of
ampleness, if $V$ is not ample, then there exists a vector subbundle $A$ of$V$ with $\deg V/A\le 0$. Since $V$ is nef, $V/A$ is nef also, and we conclude
that $\deg V/A =0$. If $A$ is the $0$ dimensional subbundle, we have $\deg V=0$, and the Lemma is proven. Indeed, if there was a nontrivial ample subbundle $A'$ in this case,
we would have $\deg V/A'=\deg V-\deg A' <0$, which contradicts the nefness of $V$. Thus we can assume that the degree of $V$ is positive. 

Thus let $A$ be a vector subbundle of $V$ of minimal rank with the property that $V/A$ has degree zero. We will be done if we show that \begin{enumerate}
\item $A$ is ample; and
\item $A$ contains any other ample subbundle of $V$. \end{enumerate} 

To see the first assertion, note that if $A$ is not ample, then by Hartshorne's theorem, there is a vector subbundle $A'$ of $A$ with $\deg A/A'\le 0$.
Thus $\deg V/A'=\deg V/A+\deg A/A'\le 0$ contradicting the choice of $A$ with minimal rank.

To see the second assertion, assume that there was an ample subbundle $B$ of $V$ with $B\not\subset A$. Then the saturation of the image of the sheaf of germs of sections of $B$ in the sheaf of germs of sections of $V/A$ is an ample subbundle of $V/A$. The inverse image of this bundle in $V$ is an ample vector subbundle of $V$ containing $A$, but the quotient of $V$ by this bundle has negative degree. This contradicts the nefness of $V$. \qed

Let $\sF$ be the ample subbundle of $\sE$ given by Lemma (\ref{hartshorne}). Clearly this bundle is also the maximal ample subbundle of $\sO_X\oplus \sE$.

Let $T:=\sE/\sF$. We claim that
$\sO_Y\oplus T^*$ is a subring of $f_*\sO_X = \sO_Y\oplus \sE^*$. Indeed, we have that $T^*$ is nef since $T$ is nef and of degree zero. Therefore, we have
that $(\sO_Y\oplus T^*)\otimes (\sO_Y\oplus T^*)$ is nef, and that $B$, the saturation of its image in $\sO_Y\oplus \sE^*$, is nef. Since $ \sO_Y\oplus T^*\subset B$,
the quotient $B/(\sO_Y\oplus T^*)$ is a nef subbundle of the bundle $\sF^*$. But since $\sF$ is ample, we conclude that $B=T^*$ and $\sO_Y \oplus T^* \subset
f_*(\O_X)$ is a subring.

Now consider the analytic spectrum $$Z := {\rm Specan}(\sO_Y \oplus T^*).$$ It comes along with a finite map
$b: Z \to Y$, and the inclusion $\sO_Y \oplus T^* \to \sO_Y \oplus \sE^*$ gives a map
$$ a: {\rm Specan}(\sO_Y \oplus \sE^*) = X \to Z $$ such that $f = b \circ a.$ Since $T^*$ equals its saturation in $\sE^*$, we conclude that $Z$ is normal and hence smooth. Since $\deg T = 0,$ we conclude that $b$ is unramified. \qed

\section{Coverings over surfaces}\label{surfaceSection} 

\begin{lemma}
Let $f: X\to Y$ be a finite morphism from an irreducible normal Cohen-Macaulay projective variety $X$, with a finite irrational locus, onto a projective manifold $Y$ of dimension $n\ge 2$. Let $H$ be a line bundle on $Y$ such that $-\omega_Y-H$ is nef and big. Then $h^1(\sE\otimes H^*)=0$. \end{lemma}

\begin{theorem}\label{delPezzoSurfaceTheorem} Let $f: X\to Y$ be a finite morphism from an irreducible
normal Gorenstein projective surface $X$ onto a smooth projective surface $Y$. Let $Y$ be a Del Pezzo surface with $\omega_Y^2 =9-r \geq 5$ (i.e., $Y$ is a smooth quadric, or the plane blown up in $r\le 4$ points in sufficiently general position). Then $\sE$ is spanned. \end{theorem}

\proof Since we already know this for the quadric, we can assume that we have a blowing up map $p : Y\to \pn 2$. We have $$ \omega_Y=-3p^*H+\sum_{i=1}^rE_i
$$
the $E_i$ are smooth $\pn 1$s on $Y$ with self-intersection $-1$. We will identify $H$ with $p^*(H).$

First let $r=1$. Let $A=2H-E_1$. Note that $-\omega_Y-A=H$ is nef and big, and thus that $h^1(\sE(-A))=0$ (3.1). We thus have $\displaystyle H^0(\sE)\to H^0(\sE_C)\to 0$. Since $A$ is very ample, a general $C \in |A-y|$ for any $y\in Y$ is smooth. For a general $C \in |A-y|$, we conclude that $f^{-1}(C)$ is reduced. Since smooth $C\in |A|$ have genus $0$, we conclude by Theorem (\ref{newLazarsfeld}) that $\sE_C$ is spanned, hence $\sE$ is spanned.

Now let $r\ge 2$. Let $H_i=H-E_i$. Note that $-\omega_Y-H_i=2H-\sum_{j\not=i}E_j$ is nef and big. Since $H_i$ is spanned with smooth $C\in |H_i|$ having genus
$0$ we conclude that $\sE$ is spanned by global sections over a generic smooth $C\in |H_i|$. Thus $\sE$ is spanned at all points $y\in Y$ where there is a smooth
$C\in |H_i-y|$ for some $i$. The set of $y\in Y$ for which this is not true for a fixed $i$ is
$B_i:=\cup_{j\not= i}\left(E_j\cup E_{ij}\right)$ where $E_{ij}$ is the proper transform under $p$ of the line through $\{p(E_i),p(E_j)\}$. Note that $\cap_i B_i=\emptyset$ unless $r=2$. 

In the case $r=2$ we have that $B_1\cap B_2=E_{12}$. Note that $-\omega_Y-H=2H-E_1-E_2$ is nef and big, and therefore we can use $H$ in place of $H_i$. This shows spannedness except at the points $y_1:=E_1\cap E_{12}$ and $y_2:=E_2\cap E_{12}$. Thus if we show that $$ H^0(\sE)\to H^0(\sE_{E_{12}})\to 0,
$$
then we will be done. To see this note that $-\omega_Y-(E_i+E_{12})$ is nef and big. Therefore using the Koszul complex and (3.1) $$
0\to \sE\otimes(-E_1-E_{12}-E_2)\to
(\sE\otimes(-E_1-E_{12}))\oplus
(\sE\otimes(-E_2-E_{12}))\to
\sE\otimes(-E_{12})\to 0,
$$
we are reduced to showing that
$h^2(\sE\otimes(-E_1-E_{12}-E_2)=0$. But we have $$ h^2(\sE\otimes(-E_1-E_{12}-E_2))=h^0(\sE^*\otimes (-2H+E_1+E_2))=h^0(f^*(-2H+E_1+E_2))=0.
$$
\qed

To get spannedness over Hirzebruch surfaces, it is enough to assume that the branch locus does not contain any components of the inverse image of $E$. \begin{theorem}\label{spannedOverHirz}Let $f:X\to Y$ be finite surjective morphism from a pure dimensional, normal variety $X$ with at worst Cohen-Macaulay singularities onto a smooth Hirzebruch surface $Y:=\hirz r$ (with surjection $ p: Y\to\pn1$). Let $E$ denote a section of $ p$ with $E\cdot E=-r$. If $r\le 1$; or $r>1$ and $E$ does not belong to the discriminant variety of $f$, then $\sE$ is spanned. \end{theorem}
\proof By Theorem (\ref{delPezzoSurfaceTheorem}), we may assume without loss of generality that $r\ge 2$.

Let $\oline E :=f^{-1}(E)$. Since $E$ does not belong to the discriminant variety of $f$, it follows that $\oline E$ is reduced. Let $d:=\deg(f)$ be the degree of $f.$

Note that global sections of $\sE$ span $\sE_E$. To see this, observe that by Theorem (\ref{newLazarsfeld}), we need only show that $h^1(\sE(-E))=0$. But this is equal to $h^1(\sE^*\otimes \omega_{Y}\otimes \sO_Y(E))$. Since $-\omega_Y-E$ is ample, we conclude that $h^1(\omega_{Y}\otimes \sO_Y(E))=0$ and $h^1(f^*(\omega_{Y}+E))=0$. Therefore
$$h^1(\sE(-E))=h^1(f^*(\omega_{Y}+ E))=0.$$ 

Given an arbitrary fiber $F$ of $ p$ we see that $\sE_{F}$ is spanned since global sections span it at $F\cap E$. Thus we will be done if we show that $h^1(\sE(-F))=0$. But this is equal to $h^1(\sE^*\otimes \omega_{Y}\otimes \sO_Y(F))$. Since $h^1(\omega_Y+F)=0$ we are reduced to showing that $h^1(f^*(\omega_Y+F))=0$. Since $f^*(\omega_Y+F)=-f^*(E+(r+1)F)-f^*E$, $f^*(E+(r+1)F)$ is ample, and $f^{-1}(E)$ is reduced, we have vanishing by Lemma (\ref{almostKodaira}). \qed

\begin{remark}
Theorem (\ref{HirzWithTwist}) is sharp for $r\ge 2$. Indeed consider the cyclic branched cover of $\hirz r$ with discriminant locus given by $E+(r-1)\left(E+rF\right)$. \end{remark} 

We have not pursued the question of what is the ``minimal'' twist, which makes $\sE$ spanned, but such questions can be answered for special classes of varieties by the same techniques. Here is a simple example. 

\begin{corollary}\label{HirzWithTwist} Let $f: X\to Y$ be a finite morphism from an irreducible normal projective surface $X$ onto smooth Hirzebruch surface $Y:=\hirz r$. If $r\ge 1$, then $\sE\otimes \sO_Y((r-1)F)$ is spanned. \end{corollary} 

\proof Let $p: Y\to\pn1$ be the $\pn 1$-bundle projection of $Y$ onto $\pn 1$.
Let $F$ be a fiber of $ p$ and let $E$ be the section of $Y$ with $E\cdot E=-r$.
To see that $\sE\otimes \sO_Y((r-1)F)$ is spanned, let $y\in Y$ be a point. A general $C\in |E+(r+1)F|$ passing through $y$ is smooth and does not belong to the branch locus. Therefore, it suffices by Theorem \ref{newLazarsfeld} to show that $H^1(\sE\otimes\sO_Y((r-1)F)\otimes\sO_Y(-C))=0$. But the latter group is isomorphic to
$H^1(f^*\sO_Y(-E-rF))$, and is therefore zero since $E+rF$ is nef and big. \qed

For higher dimensional projective bundles, it is not so easy to prove spannedness theorems under natural conditions. There is one special case for which the simple technique of splitting the anticanonical divisor works. Since we do not need this in the rest of the paper, we leave the proof to the reader.

\begin{theorem}\label{projectiveBundleOverP1} Let $f:X\to Y$ be finite surjective morphism from a pure dimensional, normal variety $X$ with finite 

irrational locus and at worst Cohen-Macaulay singularities onto $Y$, a $\pn {n-1}$-bundle over $\pn 1$.
\begin{enumerate}
\item If $-\omega_Y$ nef, then $\sE$ is spanned except if $Y = {\Bbb P}(\sO_{\pn 1}^{\oplus (n-1)} \oplus \pnsheaf 1 2),$ in which case $\sE$ is spanned outside the exceptional divisor ${\Bbb P}(\sO_{\pn 1}^{\oplus n-1}).$
\item If $Y={\Bbb P}(V)$ with
$V:=\sO_{\pn 1}\oplus \pnsheaf 1{a_1}\oplus \cdots\oplus\pnsheaf 1 {a_{n-1}}$ and with $\displaystyle \deg V\le 1+a_{n-1}$, then $\sE$ is generically spanned. \end{enumerate}
\end{theorem}

The bundle $\sE$ will fail to be ample if any curve on $Y$ has disconnected inverse image under $f$. The examples in (\ref{badExamples}) show that the following is optimal.

\begin{theorem}\label{ampleOverHirz} Let $f:X\to Y$ be a finite surjective morphism from a pure dimensional, normal variety $X$ with at worst Cohen-Macaulay singularities onto a smooth Hirzebruch surface $Y:=\hirz r$ (with surjection $ p: Y\to\pn1$). Let $E$ denote a section of $ p$ with $E\cdot E=-r$. and let $F$ denote a general fiber of $ p$. If $f^{-1}(F)$ is connected and $f^{-1}(E)$ is connected and reduced, then $H^0(\sE(-E))$ spans $\sE(-E)$ on all of $Y$ and $H^0(\sE(-F-E))$ spans $\sE(-F-E)$ on the complement of $E$. From this it follows that $\sE$ is ample. \end{theorem} 

\proof We first show that $\sE(-E)$ is spanned by global sections. Let $C$ be a general element of $|E+(r+1)F|$
that contains $y$. Since $E+(r+1)F$ is very ample, $C$ is smooth and does not belong to the discriminant variety of $f$. Thus $\oline C := f^{-1}(C)$ is reduced. It is also connected since it is an ample divisor on $X$. Thus $\sE_C(-1)$ is spanned by Theorem (\ref{newLazarsfeld}). The assertion that $\sE(-E)$ is spanned by global sections, will therefore follow if we show that $h^1(\sE(-E-C))=0$. But, the usual computation shows that this equals $h^1(f^*(-F))$. Since $f^{-1}(F)$ is connected, we conclude that the fibers of $ p\circ f:X\to \pn 1$ are all connected and thus $h^1(f^*(-F))=0$.

Now we will show that $\sE(-F-E)$ is spanned by global sections on the complement of $E$. Let $C$ be a general element of $|E+rF|$ that contains $y$. Since global sections
of $E+rF$ embed $Y-E$, the curve $C$ is smooth and does not belong to the discriminant variety. Thus $\oline C := f^{-1}(C)$ is reduced. It is also connected since it is a big divisor on $X$. Thus $\sE_C(-1)$ is spanned by Theorem (\ref{newLazarsfeld}). The assertion that $\sE(-F-E)$ is
spanned by global sections, will follow if we show that $h^1(\sE(-F-E-C))=0$. But,
the usual computation shows that this equals $h^1(f^*(-F))$. Since $f^{-1}(F)$ is connected, we conclude that the fibers of $ p\circ f:X\to \pn 1$ are all connected and thus $h^1(f^*(-F))=0$.

Since $\sE$ is spanned by Theorem (\ref{spannedOverHirz}), for $\sE$ to fail to be ample, there must exist an irreducible and reduced curve $D\subset Y$ such that
$\sE_D$ has a trivial line bundle as a direct summand. Since $\sE(-E)$ is spanned, we conclude that $D\cdot E=0$ so that $D\in |a(E+rF)|$ for some integer $a>0$. Since $\sE(-E-F)$ is spanned off $E$ and since $D\subset Y\setminus E$, we conclude that $D\cdot F=0$. Since $D\sim a(E+rF)$, we conclude that $a=0$ and thus that $D$ is the empty curve. \qed

\begin{remark} Hirzebruch \cite{Hi83} (see also, \cite{BHH87}), showed how to associate smooth projective surfaces to configurations of lines on $\pn2$. Let $\Lambda$ denote the configuration of six lines through four points in $\pn 2$, no two of which are collinear. Associated to this configuration there is a
smooth projective surfaces $\sH(\Lambda,k)$ for each integer $k\ge 2$ with the following properties:
\begin{enumerate}
\item there are finite morphisms
$f_k:\sH(\Lambda,k)\to Y$, where $Y$ is $\pn 2$ blown up at the four points; and \item the covering $f_k$ is Galois with Galois group equal to the direct sum of five copies of $\zed_k$. \end{enumerate} Since $Y$ is Del Pezzo with $\omega_Y^2=5$, Theorem (\ref{delPezzoSurfaceTheorem}) shows that the bundle $\sE$ associated to $f_k$ is spanned. This is particularly interesting because of the beautiful result of Hirzebruch that $X:=\sH(\Lambda,5)$ is a quotient of the unit ball in $\comp^2$ by a freely acting discrete group. 

It would be interesting to know in general what are the properties of the bundle $\sE$ associated to the natural finite morphisms from line configuration surfaces to blowups of $\pn 2$. In particular: \begin{question}
For which line configurations $\Lambda$ and which integers $k$ are the bundles $\sE$ spanned?
\end{question}
If there are any points where more than two lines of the configuration meet, then the surface $Y$ fibers over $\pn 1$ with the inverse image of the general fiber of this fibration disconnected in the line configuration surface, see, e.g., \cite{So84}. This shows that $\sE$ will never be ample, except possibly for the
trivial configuration where at most two lines meet in a point. \end{remark} 

\section{Coverings over Del Pezzo Manifolds}\label{delPezzoSection} 

We always let $\rho(Y)$ denote the Picard number of $Y.$ The guideline of this section is the
following

\begin{conjecture} Let $f: X \to Y$ be a finite surjective morphism of degree $d\ge 2$ between projective manifolds $X,Y$ of dimension $n$. If $Y$ is Fano with $\rho(Y) = 1,$ then $\sE$ is ample (and spanned in most cases). \end{conjecture}

First we consider the case that $Y$ is a Del Pezzo manifold of dimension $n = \dim Y \geq 3,$ i.e. $-\omega_Y = (n-1)H$ with some (ample) line bundle $H.$ Let $b := H^{n}.$ Then $1 \leq b \leq 8,$ moreover all Del Pezzo manifolds are classified (see \cite{Fu90}). In our standard situation (1.1) we obtain
\begin{eqnarray}
h^0(\sE(-H)) &=& h^0(X,{\omega_{X/Y}}\otimes f^*((n-2)H)) \label{star}\\ h^0(\sE) &=& h^0(X,{\omega_{X/Y}}\otimes f^*((n-1)H)) - 1. \label{starstar} \end{eqnarray} 

As (\ref{badExamples}) shows, $\sE$ is not always spanned if $b=1$. Though we expect spannedness for the cases $b\ge 2$, we can only prove it for $b\ge 5$.

\begin{theorem}\label{delPezzoManifoldTheorem} Let $f:X\to Y$ be a finite surjective morphism from an $n$-dimensional normal irreducible Gorenstein projective
variety $X$ onto an $n$-dimensional projective manifold $Y$. Assume that the irrational locus of $X$ is finite. Let $Y$ be a Del Pezzo manifold with $b \geq 5.$
Then $\sE$ is spanned.
\end{theorem}

\proof We proceed by induction on $n,$ noticing that $H$ is very ample. We denote a general member of $\vert H \vert $ again by $H.$ If $n = 2,$ then we conclude
by (3.2). So
suppose $n \geq 3.$ By (1.17)(4) we obtain the exact sequence $$ 0 \to H^0(\sE(-H)) \to H^0(\sE) \to H^0(\sE_H) \to 0.$$ By induction hypothesis, $\sE \vert H $ is
spanned, so it follows that $\sE$ is spanned outside a finite set $Z$. But, $H$ being very ample, through a given point we can still find a smooth member of
$\vert H \vert $, hence $Z = \emptyset.$ Notice that here we need to work in the Gorenstein category since the general $f^{-1}(H)$ need not be smooth!
\qed

In the cases $1 \leq b \leq 4,$ we do not know the spannedness of $\sE$ on the corresponding Del Pezzo surfaces, we are only able to prove much weaker results. We
start with a
general result showing that almost always $\sE$ has at least one section. 

\begin{theorem}\label{weakExistenceTheorem} Let $f:X\to Y$ be a finite surjective morphism of degree $d\ge 2$ from an $n$-dimensional normal irreducible Gorenstein
projective variety $X$ onto an $n$-dimensional projective manifold $Y$. Assume that the irrational locus of $X$ is finite. Assume that $n\ge 2$ and that $Y$
is Del Pezzo. If $bd\ge 5$, i.e., if either $b\ge 3$; or $b=2$ and $d\ge 3$; or $b=1$ and $d\ge 5$; then $h^0(\sE)\not = 0$. \end{theorem}

\proof Write $-\omega_Y = (n-1)H$ with $H$ ample. Note that if $n\ge 3$, $$
h^1(\sE(-H))=h^{n-1}(-f^*(n-2)H)=0.
$$
Since $H$ is either spanned or has at most one base point, we can reduce the proof of the result to the case when $n=2$. 

So without loss of generality we assume that $n=2$. Let $L := f^*H$. We have $$
h^0(\sE)=h^2(\sE^*\otimes \omega_Y)=h^2(-L)-1=h^0(\omega_Y+L)-1. $$ Thus we have reduced to showing that $h^0(\omega_Y+L)\ge 2$. Note that since $h^0(\omega_Y+L)=h^0(\sE)+1$, we know that $h^0(\omega_Y+L)\ge 1$. 
From this it follows
that that $\omega_X+L$ is nef, see \cite{So85}. Again using the main result of \cite{So85}, note in this case that $(\omega_X+L)\cdot L=0$ implies that $-\omega_X=L$,
which implies the absurdity that $f$ is unramified. Thus we have that: \begin{enumerate}
\item $L$ is spanned outside of a finite set; \item $\omega_X+L$ is nef and $(\omega_X+L)\cdot L\ge 2$. \end{enumerate} Noting that $L\cdot L=b\deg f$ implies
that $L\cdot L\ge 5,$
the following Lemma will complete the proof of the theorem. \begin{lemma}Let $L$ be an ample line bundle on an irreducible normal Gorenstein projective surface $X$.
Assume that $L$ is spanned outside of a finite set; $h^0(\omega_X+L)\ge 1$; and $(\omega_X+L)\cdot L\ge 2$. Then either $h^0(\omega_X+L)\ge 2$ or $L\cdot L\le 4$. \end{lemma} 

{\noindent\bf Proof of Lemma.\ } We assume that $h^0(\omega_X+L)=1$ and then argue to a contradiction.
Since $h^0(L)\ge 2$, we would have
$h^0(\omega_X+L)\ge 2$ if $h^0(\omega_X)\ge 1$. Therefore without loss of generality we can assume that
\begin{equation}\label{pgIs0}
h^0(\omega_X)=0.
\end{equation}
Since $h^1(\omega_X+L)=0$, we conclude from $$ 0\to \omega_X\to \omega_X\otimes L\to \omega_C\to 0 $$ that
\begin{equation}\label{qgRelation}
g=q+1,
\end{equation}
where $q=h^1(\sO_X) $ and $2g-2=(\omega_X+L)\cdot L$ defines the arithmetic genus $g$ of a reduced curve $C$ from $|L|$. Since $g\ge 2$ by hypothesis, we conclude further that
\begin{equation}\label{qlowerBound}q\ge 1.
\end{equation}

Assume first that the irrational locus $\sI(X)$ of $X$ is nonempty. Recall the Grauert-Riemenschneider canonical sheaf $\sK_X$ is defined as $\pi_*(\omega_{\oline X})$ where $\pi : \oline X\to X$ is a desingularization.
Following the argument of \cite{So85}, we consider the sequence: $$0\to \sK_X\otimes L\to \omega_X\otimes L \to \sS\to 0, $$ where $\sS$ is a skyscraper sheaf whose support exactly equals the set of nonrational singularities $\sS(X)$. Since $h^1(\sK_X+L)=0$, $h^0(\omega_X+L)=1$, and
the locus $\sI(X)$ is nonempty, we conclude that $h^0(\sS)=1$. Using this, equation (\ref{pgIs0}), and the sequence: $$0\to \sK_X\to \omega_X \to \sS\to 0,
$$
we conclude that:
\begin{equation}\label{onDeSing}
h^0(\omega_{\oline X})=0;\ \ \ \ \
\oline q := h^1(\sO_{\oline X})=q+1.
\end{equation}
Therefore $X$ has only one irrational singularity, an elliptic singularity $y_0.$ Consider the Albanese mapping $\alpha : \oline X\to Z$. By equation (\ref{onDeSing})
we conclude first that $\oline q\ge 1$ and therefore that $\dim Z\ge 1$. Since $h^0(\omega_{\oline X})=0$, we conclude that $Z$ is a curve, and thus smooth with $\alpha$ having connected fibers. Since $q\ge 1$ by equation (\ref{qlowerBound}), we
conclude
that $\oline q\ge 2$. Let $B$ denote the singular set of $X$, and note that $\alpha(\pi^{-1}(B))$ is finite.
This is clear since $y_0$ is elliptic and all other singularities are rational. In particular, the strict transform $\overline C$ of $C$ cannot be a
fiber of $\alpha.$ Thus a general fiber $F$ of $\alpha$ is mapped isomorphically by $\pi$ onto a curve $F'$ contained in the smooth points of $X$.
Since $L$ is spanned off of a finite set, we conclude that global sections of $L$ give rise to at least two linearly independent sections of $L_{F'}$. Hence we deduce easily that $C \cdot F' \geq 2$ and $\overline C$ has degree at least 2 over $Z.$ Thus using that ${\overline {q}}\ge 2$, we have that
$$ g\ge 2\oline q \ge 2q+2=2g,$$
which is absurd since $g\ge 2$.

Therefore we can assume without loss of generality that the singularities of $X$ are rational. Thus, e.g., \cite{BS95}, the Albanese map $\alpha : X\to Z$ is well
defined. Since $h^0(\omega_X)=0$ by equation (\ref{pgIs0}), we conclude that by the same argument as above, but without the need of a desingularization, that if $q\ge 2$ then $$ g\ge 2 q =2g-2,$$
which is absurd given that $g=q+1\ge 3$. Thus we conclude that $q=1$ and $g=2$. In this case we have that
$$ (\omega_X+L)^2 L^2\le [L\cdot (\omega_X+L)]^2=1. $$ Since $L\cdot L\ge 5$ we conclude using the nefness of $\omega_X+L$ that $(\omega_X+L)^2=0$.
This gives $\omega_X^2\ge 1$.
Since $q=1$, $h^0(\omega_X)=0$, and
we have at worst rational Gorenstein singularities, we conclude from the Riemann-Roch formula for the Euler characteristic of the structure sheaf of the desingularization,
that $\omega_X^2\le 0$.

This proves the Lemma and the Theorem. \qed 

\begin{remark}Note if $Y$ is an elliptic curve, which can be considered as the most appropriate definition of a Del Pezzo manifold of dimension $1$, then Theorem
(\ref{delPezzoManifoldTheorem})
is still true, unless $f$ is an unramified covering. \end{remark} 

We now show how adjunction theory can be used to give general conditions for an existence theorem, which would not follow even from knowing that $\sE$ is spanned.
Indeed, as shown in (\ref{exceptions}), some of the possible exceptions to this result exist and in fact have $\sE$ spanned. We do not need
the precise classification of Del Pezzo manifolds with $1 \leq b \leq 4$ after Fujita \cite{Fu90}, but just the following 

\begin{proposition}\label{classification} Let $Y$ be a Del Pezzo manifold of dimension at least
3. Then $\rho(Y) = 1,$ unless $X$ is $\left(\pn 1\right)^3$, $\pn2 \times \pn 2$, $\bP(T_{\pn 2})$ or $\pn 3$ blown up in one point. Moreover $H $ is very ample for $b\ge 3$, and spanned unless $b = 1$, in which case we have
one simple base point. In all cases $\vert H \vert $ contains a smooth element. \end{proposition}

\begin{theorem}\label{minuHExistenceTheorem} Let $f:X\to Y$ be a finite surjective morphism of degree $d\ge 2$ from an $n$-dimensional projective manifold $X$ onto an $n$-dimensional projective manifold $Y$. Assume that $n\ge 3$ and that $Y$
is Del Pezzo, i.e., $-\omega_Y=(n-1)H$ with $H$ ample. Let $Y$ be of degree $b:= H^n\ge 2$.
Then $H^{0}(\sE(-H)) \ne 0$, except for the following possible cases: in the cases 3) and 4), $(X_0,L_0)$ denotes the first reduction of $(X,L)$ \cite{BS95}.
\begin{enumerate}
\item $(X,L)$ is a quadric fibration over a curve with $Y=\left(\pn 1\right)^3$; or
\item $(X,L)$ is a scroll over a surface with either $n = 3$; or \item $(X_0,L_0)$ is a 3-dimensional quadric, $L_0 = \sO(2)$, $Y =\pn3(x_0)$, $d=2$, and
$f$ is induced by a covering $f_0: X_0 \to Y_0 = \pn3$; or \item $n = 3$ and $X_0$ is a $\pn 2$-bundle over a curve with $L_F = \pnsheaf 2 2$ on a fiber $F$. \end{enumerate} Examples of dimension $3,4$, for the cases of the Theorem are given in (\ref{exceptions}).
\end{theorem}

\proof Let $L = f^*(H).$ Suppose $H^0(\sE(-H)) = 0.$ Then by (\ref{star}) $$ H^0(\omega_X \otimes f^*((n-2)H)) = 0.$$ Then, since $b\ge 2$ and $H$ is spanned, \cite{So89} applies.
Using the equation
(\ref{starstar}) that $h^0(\omega_X+(n-1)L)=h^0(\sE)+1\ge 1$, we obtain that either $(X,L)$ is
\begin{enumerate}
\item[(a)] a quadric fibration over a curve; or \item[(b)] a scroll over a surface; or,\\ \ \ \\
$(X_0,L_0)$, the first reduction
of $(X,L)$ \cite{BS95}, exists and is either: \item[(c)] $(X_0,L_0) = \pnpair n2$ with
$n = 3,4$; or
\item[(d)] $(X_0,L_0)$ is a 3-dimensional quadric and $L_0 =\sO(2);$ or \item[(e)] $n = 3$ and $X_0$ is a $\pn 2$-bundle over a curve with $L = \sO(2)$ on the fibers.
\end{enumerate}

We will treat these cases separately. \\ 

(a) If we have a quadric fibration $p: X \to C,$ let $F$ be a general fiber and $F' = f(F).$ Since $L^{n-1} \cdot F = 2,$ the map $f_F$ has degree 2 or 1.
$H$ being ample and spanned (recall $b \geq 2$), it follows that in the first case $F'$ must be projective space, and the adjunction formula implies
$N_{F' \vert Y} = \sO(-1),$ which is absurd since $F'$ moves in $Y.$ In the second case, we obtain $H^{n-1} \cdot F' = 2$ and we conclude similarly that $F'$ is a quadric. By adjunction, $N_{F' \vert Y} = \sO,$ hence $\rho(Y) \geq 2.$
By the classification (\ref{classification}) we obtain $Y = \left(\pn 1\right)^3.$ This is in one of the exceptions. \\ 

(b) Now suppose that $p: X \to S$ is a scroll over the surface $S.$ Then, using the same notations as before, we conclude $F' = \pn {n-2}$ and $H \vert F' = \sO(1).$ By adjunction, $ \det N_{F'} = \sO.$ Notice that $N_{F'}$ is generically generated for {\em general} $F'$, hence $N_{F'}$ must be trivial.
If $b \geq 3$ and $n \geq 5$, then $H$ is very ample, so that by Ein \cite{Ei85}, $Y$ has a $\pn {n-2}$-bundle structure, contradicting (\ref{classification}) unless
$Y = \pn 2 \times \pn 2 $ giving rise to one of the exceptions (4.10). Therefore it remains to treat the case $b = 2$ or $n \le 4.$ \\ (b.1) In case $n = 3,$ then $F'$ is a line and we are in the exceptions (4.10), the cases $Y = \bP(T_{\bP_2}), \bP_3(x), \bP_1^3$ being obvious.\\ (b.2) Concerning $b = 2,$ the linear system $\vert H \vert$ realizes $Y$ as a $2:1-$covering $h: Y \to \bP_n,$ moreover $L = f^*(H).$ Since $H \vert F' = \sO(1),$ it follows that $f \vert F$ is $1:1$ for all $F$. This is clearly absurd, since $f^{-1}(F')$ would always consist of 1 or 2
disjoint fibers of $p$, which implies $\rho(Y) \geq 2,$ a contradiction to the classification. \\
(b.3) It remains to treat the case $b \geq 3$ and $n = 4.$ \\ ($\alpha$) Suppose $b = 3.$ Then $Y \subset \bP_5$ has degree and we obtain (for general $F' = \bP_2$) an exact sequence $$
0 \to N_{\bP_2 \vert Y} \to N_{\bP_2\vert \bP_5} \to N_{Y \vert \bP_5} \vert \bP_2 \to 0,
$$
which reads
$$ 0 \to \sO \oplus \sO \to \sO(1)^{\oplus 3} \to \sO(3) \to 0. $$ This sequence is of course absurd. \\
($\beta$) Suppose $b = 4.$ Then $Y$ is the intersection of 2 quadrics in $\bP_6.$
Take one of the quadrics, say $Q$, to obtain a sequence $$
0 \to N_{\bP_2 \vert Y} \to N_{\bP_2 \vert Q} \to N_{Y \vert Q} \vert \bP_2 \to 0.
$$
This sequence immediately implies
$$ N_{\bP_2 \vert Q} = \sO\oplus \sO \oplus \sO(2), $$ which clearly contradicts the 1-ampleness of the tangent bundle $T_Q.$ \\ ($\gamma$) Finally we treat $b = 5.$ First notice that every line $l \subset Y $ has a
normal bundle either of the form $\sO \oplus \sO \oplus \sO(1)$ or $\sO(1) \oplus \sO(1) \oplus \sO(-1).$ This follows by choosing a smooth hyperplane 
section through $l$ (which is easily seen to exist, cp. [Is77,p.505]). We are going to show that {\it every} plane $F' \subset Y$ has normal bundle $\sO \oplus \sO.$ Once we know this, we can argue as in [Ei85,1.7] to obtain a $\bP_2-$bundle structure on $Y.$ Since the planes $F'$ cover $Y,$ the
normal bundle $N$ is generically generated by 2 sections $s_1,s_2.$ Since by adjunction $\det N = \sO,$ it follows that these sections generate $N$ everywhere so that $N$ is trivial.

(c) Suppose $(X_0,L_0) = \pnpair n 2$ with $n = 3 $ or $4.$ We claim that always $X = X_0.$
In fact, if $X \not = X_0,$ then $X$ contains a divisor $E = \pn {n-1}$ with $N_E = \sO(-1).$ Consider $E' = f(E).$ Since $L \vert E = \sO(1),$ we have $H^{n-1} \cdot E' = 1,$ so that since $H$ is therefore spanned, we have that $E' = \pn {n-1}.$ Adjunction gives $N_{E'} = \sO(-1),$ hence by classification $Y = \pn 3(x_0).$ Therefore $X_0 = \pn 3$.
Now, introducing the ramification divisor $R$ of $f,$ we have $$\omega_X + f^*(H) = f^*(\omega_Y) + R + f^*(H) = f^*(-H) + R,$$ hence $$ H^0(X,\sO_X(R) \otimes f^*(-H)) = 0$$ by our assumption. If $\sigma: X \to X_0$ denotes the contraction of the exceptional divisors $E,$ we conclude that $$R' = \sigma_*(R) \in \vert \sO_{X_0}(1) \vert.$$ Hence the induced covering $h: X_0 \to Y_0$ is ramified exactly over a plane, which is absurd since $X_0 = Y_0 = \pn 3.$ So we shall assume $X = X_0.$ But then we have the contradiction that $Y = \pn n$.

(d) Suppose $X_0$ is a three dimensional quadric with $L=\sO_{X_0}(2)$. If $X=X_0$, then $Y$ is $\pn 3$ or a quadric. Since $Y$ is Del Pezzo, we have
$Y=\pn 3$ and $H=\pnsheaf 3 2$. Since $K_X+L$ has no sections, we conclude that $L=\sO_{X}(2)$, and therefore that $f$ has degree $2$. From this we conclude
that $\sE=\pnsheaf 3 1$. This gives an example. If $X \not = X_0,$ then we argue as in (c) to conclude that $Y = \pn(3)(x_0).$ Hence $Y_0 = \pn3,$ see (4.10)(4).
\qed

\begin{remark}In the preceding and following results we need $b\ge 2$ to get spannedness so that the main result of \cite{So89} holds. Actually, we only need that $2L$ is spanned and $|L|$ has a smooth element. This would be satisfied when $b=1$ under the extra condition that the point where $H$ is not spanned is not an element of the discriminant divisor of the covering . Using this the statement of the result is the same, except there is one case that cannot be ruled out: the first reduction $(X_0,L_0)$ of $(X,L)$ exists and is equal
$\pnpair n2$ with
$n = 3,4$.
\end{remark}

\begin{corollary} If $b\ge 2$, then $h^0(\sE) \geq h^0(\sO_Y(H)) $ unless $X$ is of the one of the exceptions in Theorem (\ref{minuHExistenceTheorem}). In particular these inequalities hold if $n \geq 5.$
\end{corollary}

We now give examples of the exceptional cases. 

\begin{exceptions}\label{exceptions} {\rm (1) $Y = \left(\pn 1\right)^3$ and we fix a projection $q:
\left(\pn 1\right)^3 \to \pn 1.$ Let $h: C \to \pn 1$ be a suitable ramified cover and put $X = Y \times_{ \pn 1} C.$ \\ 

(2) Let $Y$ be a Del Pezzo 3-fold and $T$ a component of the Hilbert schemes of lines. Suppose that through every point of $X$ there are only finitely many lines. Let
$S$ be a smooth surface and $S \to T$ be a surjective map. Let $X$ be the $ \pn 1$-bundle over $S$ given by the graph of the family of line parameterized by $T.$ \\

(3) Let $Y = \pn 2 \times \pn 2.$ Let $S$ be a surface and $S \to \pn 2$ be a finite map.
Put $X = S \times \pn 2 \to Y.$ Analogously we choose $S \to \pn 1 \times \pn 1$ finite and let $X = S \times \pn 1$ resp. $h: S \to \pn 2$ finite and $X = \bP(h^*(T_{ \pn 2})).$ \\

(4) Suppose $X_0$ is a three dimensional quadric. We have a degree $2$ map $f: X_0\to \pn 3$. In this case $\sE=\pnsheaf 3 1$. Since $\pn 3$ is Del Pezzo
with $H=\pnsheaf 3 2$, we have $\sE(-H)$ has no sections. Note $L_0=\sO_{X_0}(2)$.
If $\pi : Y\to \pn 3$ is the blowup of $\pn 3$ at a point $y$. Then $\oline H:=\pi^*\pnsheaf 3 2-E$ is very ample where $E$ is the exceptional fiber of $\pi$.
Thus $-K_Y=\pi^*\pnsheaf 3 4-2E=2\oline H$, and we obtain a Del Pezzo manifold $Y$. Letting $X$ denote the blowup of $X$ at the inverse image under $f$ of $y$,
we have a degree $2$ finite morphism $\oline f : X\to Y$. We have $h^0(K_X+L)=h^0(K_{X_0}+L_0)=0$. \\

(5) Let $p: X_0 \to C$ be a $ \pn 2$-bundle over a smooth curve $C.$ Let $h: X_0 \to Y_0 = \pn 3$ be a finite
cover such that $$[h^*{\pnsheaf 3 1}]_F = \pnsheaf 2 1,$$ where $F$ denotes a fiber of $p.$ Choose $y_0$ not in the image of the ramification divisor of $h.$
Then let $\tau : Y \to Y_0$ be the blow-up of $Y_0$ at $y_0$ and $\sigma: X \to X_0$ the
blow-up of $X_0$ at the finite set $h^{-1}(y_0).$ Then $h$ lifts to a finite cover $f: X \to Y$ and it is easy to check that $H^0(X,\omega_X \otimes f^*(H)) = 0,$
where $-\omega_Y = 2H.$

In order to obtain a specific example, consider the family ${\cal C}$ of hyperplanes in $\pn 3.$ The canonical
map $$q: {\cal C} \to \pn 3 = Y_0$$ is a $\pn 2$- bundle and so does the other projection $$p: {\cal C} \to {\pn 3}^*.$$ Actually ${\cal C} = \bP(T_{\pn 3}).$ Now
choose $C \subset {\pn 3}^*$ to be a general smooth curve. Then $X_0 = p^{-1}(C)$ is smooth and moreover $h := q \vert X_0$ is finite. To see this last property,
consider $T(x) \subset {\pn 3}^*,$ the set of hyperplanes through a given point $x \in \pn 3$. So it is sufficient to choose $C$ such that $C \not \subset T(x)$ for
all $x$ which comes down to choose $C$ non-degenerate, since all $T(x)$ are linear in $({\pn 3})^*.$ }
\end{exceptions}

\vspace{1cm}
\small
\begin{tabular}{lcl}
Thomas Peternell&& Andrew J. Sommese \\ Mathematisches Institut & & Department of Mathematics\\ Universit\" at Bayreuth & & University of Notre Dame\\ D-95440 Bayreuth, Germany&&Notre Dame, Indiana 46556, U.S,A,\\ fax: Germany + 921--552999&& fax: U.S.A. + 219--631-6579 \\ thomas.peternell@uni-bayreuth.de && sommese@nd.edu\\ &&URL: {\tt www.nd.edu/$\sim$sommese}\\
\end{tabular}


\begin{thebibliography}{999999}

\bibitem[BHH87]{BHH87} Barthel,G.; Hirzebruch,F.; H\"ofer,Th.: Geradenkonfigurationen und Algebraische Fl\"achen. Aspects of Mathematics, D4. Friedr. Vieweg \& Sohn 1987


\bibitem[BS95]{BS95}Beltrametti,M.; Sommese,A.J.: The adjunction theory of complex projective varieties. de Gruyter 1995 

\bibitem[CP97]{CP97}Campana,F.;Peternell,Th.: Towards a Mori theory on compact K\"ahler threefolds, I. Math. Nachr. 187, 29--59 (1997) 


\bibitem[CS94]{CS94}Cho,K.;Sato,E.: Smooth projective varieties dominated by smooth hyperquadrics in any characteristic. Math. Z. 217, 553--565 (1994) 


\bibitem[De95]{De95}Debarre,O.: Th\'eoremes de connexit\'e et vari\'et\'es ab\'eliennes. Amer.J.Math. 117, 787--805 (1995) 

\bibitem[De96]{De96} Debarre,O.: Sur un th\'eoreme de connexit\'e de Mumford pour les espace homogenes. Manuscr.\ Math. 89, 407-425 (1996) 

\bibitem[DPS94]{DPS94} Demailly,J.P.;Peternell,Th.;Schneider,M.: Compact complex manifolds with numerically effective tangent bundles. J. Alg. Geom. 3, 295--345 (1994)

\bibitem[DPS99]{DPS99} Demailly,J.P.;Peternell,Th.;Schneider,M.: Pseudo-effective line bundles on projective varieties and K\"ahler manifolds. Preprint in prep. (1999)

\bibitem[Ei85]{Ei85} Ein,L.: Varieties with small dual varieties,II. Duke Math. J. 52, 895--907 (1985)

\bibitem[Fu90]{Fu90} Fujita,T.: Classification theories of polarized varieties. London Math. Soc. Lect. Notes Ser. 155 (1990) 

\bibitem[Ha71]{Ha71} Hartshorne,R.: Ample vector bundles on curves. Nagoya Math.\ J. 43, 73--89 (1971)

\bibitem[Hi83]{Hi83} Hirzebruch,F.: Arrangements of lines and algebraic surfaces. Arithmetic and geometry, Vol. II, 113--140, Progr. Math., 36, Birkh\"auser, Boston, Mass., 1983


\bibitem[Is77]{Is77} Iskovskih,V.A.: Fano 3-folds I. Math. USSR Izv. 11, 485-527 (1977)

\bibitem[Ki96a]{Ki96a} Kim,M.: Barth-Lefschetz type theorem for branched coverings of Grassmannians, J. Reine Angew.\ Math. 470, 109--122 (1996) 

\bibitem[Ki96b]{Ki96b} Kim,M.: On branched coverings of quadrics. Arch. Math. 67, 76--79 (1996)

\bibitem[KM98]{KM98} Kim,M.;Manivel,L.: On branched coverings of some homogeneous spaces. Topology 38, 1141--1160 (1999) 

\bibitem[Ma97]{Ma97} Manivel,L.: Vanishing theorems for ample vector bundles. Invent.\ Math.\ 127, 401--416 (1997) 


\bibitem[La80]{La80} Lazarsfeld,R.: A Barth type theorem for branched coverings of projective spaces. Math. Ann. 249, 153--162 (1980) 

\bibitem[OSS80]{OSS80} Okonek,C.;Schneider,M.;Spindler,H.: Vector bundles on complex projective spaces. Birkh\"auser 1980 

\bibitem[PS89]{PS89} Paranjap\'e,K.H.;Srinivas,V.: Selfmaps of homogeneous spaces. Inv.\ Math. 98, 425--444 (1989)

\bibitem[So78]{So78} Sommese,A.J.: Submanifolds of Abelian varieties. Math.\ Ann.\ 233, 229--256 (1978)



\bibitem[So84]{So84} Sommese,A.J.: On the density of ratios of Chern numbers of algebraic surfaces. Math.\ Ann.\ 268, 207--221 (1984) 

\bibitem[So85]{So85} Sommese,A.J.: Ample divisors on normal Gorenstein surfaces. Abh. Math. Sem. Univ. Hamburg 55, 151--170 (1985) 

\bibitem[So86]{So86} Sommese,A.J.: On the adjunction theoretic structure of projective varieties.
In: Proceedings of the Complex Analysis and Algebraic Geometry Conference, G\"ottingen, 1985, ed.\ by H. Grauert,
Lecture Notes in Math. 1194, 175--213. Springer 1986 

\bibitem[So89]{So89} Sommese,A.J.: On the nonemptiness of the adjoint linear system of
a hyperplane section of a threefold. J. reine u. angew.\ Math.\ 402, 211--220 (1989); Erratum, 411, 122--123 (1990) 


\bibitem[Vi82]{Vi82} Viehweg,E.: Die Additivit\"at der Kodairadimension f\"ur projektive Faserr\"aume \"uber Varie\"aten des allgemeinen Typs. J.reine u. angew. Math. 330, 132--142 (1982) 

\bibitem[Vi95]{Vi95} Viehweg,E.: Quasi-projective moduli for polarized manifolds. Springer 1995

\end{thebibliography}
\end{document}